\documentclass[11pt]{article} \makeatother
\usepackage{amsfonts,amsmath,amssymb}

\begin{document}

\numberwithin{equation}{section}
\newtheorem{theorem}{\ \ \ \ Theorem}[section]
\newtheorem{proposition}[theorem]{\ \ \ \ Proposition}
\newtheorem{lemma}[theorem]{\ \ \ \ Lemma}
\newtheorem{remark}{\ \ \ \ Remark}
\newcommand{\be}{\begin{equation}}
\newcommand{\ee}{\end{equation}}
\newcommand\bes{\begin{eqnarray}} \newcommand\ees{\end{eqnarray}}
\newcommand{\bess}{\begin{eqnarray*}}
\newcommand{\eess}{\end{eqnarray*}}
\newcommand\ds{\displaystyle}

 \begin{center} {\bf\Large Exponential energy decay of solutions for a system of viscoelastic}\\[2mm]
  {\bf\Large     wave equations of Kirchhoff type with strong damping}
\\[4mm]
  {\large Gang Li, Linghui Hong, Wenjun Liu$^{1}$}\\[1mm]
{\small  College of Mathematics and Physics, Nanjing University of
Information Science and Technology, Nanjing 210044, China. E-mail: wjliu@nuist.edu.cn.  }\\[1mm]
\end{center}

\setlength{\baselineskip}{17pt}{\setlength\arraycolsep{2pt}

\begin{quote}
\noindent {\bf Abstract:} The initial boundary value problem for a
system of viscoelastic wave equations of Kirchhoff type with strong
damping is considered. We prove that, under suitable
 assumptions on relaxation functions and certain
 initial data, the decay rate of the solutions energy is
 exponential.


\end{quote}

\setlength{\baselineskip}{17pt}{\setlength\arraycolsep{2pt}

\section{Introduction}

In this work, we investigate the following system of viscoelastic
wave equations of Kirchhoff type:
 \bes\left\{\begin{array}{ll}
\displaystyle  u_{tt}-M(\|\nabla u\|_{2}^{2})\Delta
u+\int_{0}^{t}g_{1}(t-s)\Delta u(s){\rm d}s-\Delta
u_{t}=f_{1}(u,v),\quad &(x,t)\in \Omega \times (0,\infty
),\medskip\\\medskip
 \displaystyle  v_{tt}-M(\|\nabla v\|_{2}^{2})\Delta v+\int_{0}^{t}g_{2}(t-s)\Delta
v(s){\rm d}s-\Delta v_{t}=f_{2}(u,v), &(x,t)\in \Omega \times
(0,\infty ), \\\medskip
  \displaystyle u(x,t)=0, \quad  v(x,t)=0, &(x,t)\in \partial\Omega \times [0,\infty ),\\\medskip
\displaystyle u(x,0)=u_{0}(x),\quad u_{t}(x,0)=u_{1}(x), \quad   &x
\in  \Omega, \\\medskip
 \displaystyle v(x,0)=v_{0}(x),\quad v_{t}(x,0)=v_{1}(x),
\quad   &x \in  \Omega,
 \end{array}\right.\label{1.1}
 \ees
where $\Omega\subset  \mathbb{R}^{n},$ $n\geq 1,$ is a bounded
domain with smooth boundary $\partial\Omega,$ $M$ is a positive
$C^{1}$ function and $g_{i}(\cdot):\mathbb{R}^{+}\rightarrow
\mathbb{R}^{+},$ $f_{i}(\cdot,\cdot):\mathbb{R}^{2}\rightarrow
\mathbb{R}^{2} (i=1,2)$ are given functions to be specified later.

The motivation of our work is due to some results regarding
viscoelastic wave equations of Kirchhoff type. The single wave
equation of the form
\begin{equation}
u_{tt}-M(\|\nabla u\|_{2}^{2})\Delta u+\int_{0}^{t}g(t-s)\Delta
u(s){\rm d}s+h(u_{t})=f(u),\quad (x,t)\in \Omega \times (0,\infty ),
\label{1.2}
\end{equation}
is a model to describe the motion of deformable solids as hereditary
effect is incorporated. Equation \eqref{1.2} was first studied
by Torrej\'on and Young \cite{21a} who proved the existence of
weakly asymptotic stable solution for large analytical datum. Later,
Rivera \cite{14a} showed the existence of global solutions for small
datum and the total energy decays to zero exponentially under some
restrictions. Then, Wu and Tsai \cite{2006a} treated equation \eqref{1.2}
for $h(u_{t})=-\Delta u_{t},$ they established the global existence
as well as energy decay under assumption on the nonnegative kernel
$g'(t)\leq -rg(t), \forall t\geq0$ for some $r>0.$ This energy decay
result was recently improved by Wu in \cite{2009a} under a weaker
condition on $g$ (i.e., $ g'(t)\leq0$ for $t\geq0$). For a single
wave equation of Kirchhoff type that without the viscoelastic term, we refer the reader to Refs.
\cite{18c, 14c, 15c, 16c, 17c}.

The system of wave equations that without viscoelastic terms $(i.e.,
g_{i}=0, i=1,2)$ has also been extensively studied and many results
concerning local, global existence, decay and blow-up have been
established. For example, Park and Bae \cite{23e} considered the
system of wave equations with nonlinear dampings for
$f_{1}(u,v)=\mu|u|^{q-1}u$ and $f_{2}(u,v)=\mu|v|^{q-1}v,$ $q\geq1,$
and showed the global existence and asymptotic behavior of solutions
under some restriction on the initial energy. Later, Benaissa and
Messaoudi \cite{3e} discussed blow-up properties for negative
initial energy. Recently, Wu and Tsai \cite{2007a} studied the
system \eqref{1.1} for $g_{i}=0\ (i=1,2).$ Under some suitable
assumptions on $f_{i}\ (i=1,2),$ they proved the local existence of
solutions by Banach fixed point theorem and blow-up of solutions by
using the  method of Li and Tsai in \cite{2003a}, where three
different cases on the sign of the initial energy $E(0)$ are
considered.

For the case of $M\equiv 1$ and in the presence of the viscoelastic
terms $(i.e., g_{i}\neq 0, i=1,2),$   problem
\bes\left\{\begin{array}{ll} \displaystyle u_{tt}-\Delta
u+\int_{0}^{t}g_{1}(t-\tau)\Delta u(\tau){\rm d}\tau-\Delta
u_{t}=f_{1}(u,v),\quad &(x,t)\in \Omega \times (0,T
),\medskip\\\medskip
 \displaystyle  v_{tt}-\Delta v+\int_{0}^{t}g_{2}(t-\tau)\Delta
v(\tau){\rm d}\tau-\Delta v_{t}=f_{2}(u,v), &(x,t)\in \Omega \times
(0,T),
\\\medskip
  \displaystyle u(x,t)=0, \quad  v(x,t)=0, &(x,t)\in \partial\Omega \times [0,T ),\\\medskip
\displaystyle u(x,0)=u_{0}(x),\quad u_{t}(x,0)=u_{1}(x), \quad   &x
\in  \Omega, \\\medskip
 \displaystyle v(x,0)=v_{0}(x),\quad v_{t}(x,0)=v_{1}(x),
\quad   &x \in  \Omega,
 \end{array}\right.\label{1.3}
 \ees
 was recently studied by Liang and Gao in \cite{2011a}.
Under suitable assumptions on the functions $g_{i},$ $f_{i}$
$(i=1,2)$ and certain initial data in the stable set, they proved
that the decay rate of the solution energy is exponential.
Conversely, for certain initial data in the unstable set, there are
solutions with positive initial energy that blow up in finite time.

It is also worth mentioning the work \cite{2009b} in which the authors considered the following problem
\bes\left\{\begin{array}{ll} \displaystyle u_{tt}-\Delta
u+\int_{0}^{t}g_{1}(t-\tau)\Delta u(\tau){\rm d}\tau+
|u_{t}|^{m-1}u_{t}=f_{1}(u,v),\quad &(x,t)\in \Omega \times (0,T
),\medskip\\\medskip
 \displaystyle  v_{tt}-\Delta v+\int_{0}^{t}g_{2}(t-\tau)\Delta
v(\tau){\rm d}\tau+|v_{t}|^{\gamma-1}v_{t}=f_{2}(u,v), &(x,t)\in
\Omega \times (0,T), \\\medskip
  \displaystyle u(x,t)=0, \quad  v(x,t)=0, &(x,t)\in \partial\Omega \times [0,T ),\\\medskip
\displaystyle u(x,0)=u_{0}(x),\quad u_{t}(x,0)=u_{1}(x), \quad   &x
\in  \Omega, \\\medskip
 \displaystyle v(x,0)=v_{0}(x),\quad v_{t}(x,0)=v_{1}(x),
\quad   &x \in  \Omega,
 \end{array}\right.\label{1.4}
 \ees
 where $\Omega$ is a bounded domain with smooth boundary
 $\partial\Omega$ in $\mathbb{R}^{n}, n=1,2,3.$ Under suitable assumptions on
 the functions $g_{i},$ $f_{i}$ $(i=1,2),$ the initial data and the
 parameters in the above problem, they established local existence,
 global existence and blow-up
 property (the initial energy $E(0)<0$). This latter blow-up result has been improved by Messaoudi
 and Said-Houari \cite{2010b}, into certain solutions with positive
 initial energy.
For other papers related to existence, uniform decay and blow-up of
solutions of nonlinear wave equations, see \cite{5e, 11c, 13e, 8e,
lna, ljmp, 9c, 24e, 7e, sw2006, 16b, 25e, 2009a, y2009} and
references therein.

Motivated by the above mentioned research, we consider in the
present work the coupled system \eqref{1.1} with nonzero $g_{i}$ $(i=1,2)$ and nonconstant $M(s)$. We show that,
under suitable assumptions on the functions $g_{i},$ $f_{i}$
$(i=1,2)$ and certain initial conditions, the solutions are global
in time and the energy decays exponentially.

This paper is organized as follows. In section 2, we first give some
assumptions, notations and lemmas which will be used later, and then
state the local existence. In section 3, we
devote to state and prove our main result.

\section{Preliminaries and main result}

In this section we first present some assumptions, notations and lemmas
needed for our work, and then state the local existence theorem.
First, we make the following assumptions:
 \begin{quote}
(A1) $M(s)$ is a positive $C^{1}$ function for $s\geq0$ satisfying
$$M(s)=m_{0}+m_{1}s^{\gamma},\quad m_{0}>0,\quad m_{1}\geq0\quad \text{and} \quad  \gamma\geq1.$$
(A2) $g_{i}(t):\mathbb{R}^{+}\rightarrow \mathbb{R}^{+}$ $(i=1,2)$ belong
to $C^{1}(\mathbb{R}^{+})$ and satisfy $$g_{i}(t)\geq 0, \quad
g_{i}'(t)\leq 0, \quad \text{for}\quad t\geq 0.$$
\end{quote}

%

Next, we introduce some notations:
\begin{align*}
k_{i}=m_{0}-\int_{0}^{\infty}g_{i}(s){\rm d}s>0, \quad
k=\min\{k_{1},k_{2}\},
\end{align*}
\begin{align*}
(g_{i} \circ \nabla w)(t) =\int_{0}^{t}g_{i} (t-s)\|\nabla
w(t)-\nabla w(s)\|_{2}^{2}\,{\rm
  d}s,
\end{align*}
\begin{align*}
g_{0}=\max\{g_{1},g_{2}\}, \quad
\|\cdot\|_q=\|\cdot\|_{L^q(\Omega)}, \quad 1\le q\le \infty,
\end{align*} the Hilbert space $L^{2}(\Omega)$ endowed with the inner
product
$$(u,v)=\int_{\Omega}u(x)v(x){\rm d}x,$$
and the functions $f_{1}(u,v)$ and $f_{2}(u,v)$ (see also \cite{2010b})
$$f_{1}(u,v)=\left[a|u+v|^{2(p+1)}(u+v)+b|u|^{p}u|v|^{p+2}\right],$$
$$f_{2}(u,v)=\left[a|u+v|^{2(p+1)}(u+v)+b|u|^{p+2}|v|^{p}v\right],$$
where $a, b>0$ are constants and $p$ satisfies
 \bes\left\{\begin{array}{ll}
\displaystyle   -1<p, \quad &\text{if}\quad n=1,2,
\medskip\\\medskip
 \displaystyle  -1<p\leq (3-n)/(n-2), \quad &\text{if}\quad n\geq3.
\end{array}\right. \label{2.1}
 \ees

 One can easily verify that
$$uf_{1}(u,v)+vf_{2}(u,v)=2(p+2)F(u,v), \quad \forall\ (u,v)\in \mathbb{R}^{2},$$
where
$$F(u,v)=\frac{1}{2(p+2)}\left[a|u+v|^{2(p+2)}+2b|uv|^{p+2}\right].$$

The energy functional $E(t)$ and auxiliary functional $I(t),$ $J(t)$
of the solutions $(u(t),v(t))$ of problem \eqref{1.1} are defined as
follows
\begin{align}
 I(t):=&I(u(t),v(t))=\left(m_{0}-\int_{0}^{t}g_{1}(s){\rm
d}s\right)\|\nabla u(t)\|_{2}^{2}
+\left(m_{0}-\int_{0}^{t}g_{2}(s){\rm
d}s\right)\|\nabla v(t)\|_{2}^{2}\nonumber\\
&+(g_{1}\circ \nabla u)(t)+ (g_{2}\circ \nabla
v)(t)-2(p+2)\int_{\Omega}F(u,v){\rm d}x,\label{2.2}
\end{align}
\begin{align}
 J(t):=&J(u(t),v(t))=\frac{1}{2}\left[\left(m_{0}-\int_{0}^{t}g_{1}(s)\, {\rm
d}s\right)\|\nabla u\|_{2}^{2} +\left(m_{0}-\int_{0}^{t}g_{2}(s)\,
{\rm d}s\right)\|\nabla v\|_{2}^{2}\right]-\int_{\Omega}F(u,v){\rm
d}x
\nonumber\\
&+\frac{1}{2}\left[(g_{1}\circ \nabla u)(t)+ (g_{2}\circ \nabla
v)(t) +\frac{m_{1}}{\gamma+1}\left(\|\nabla u\|_{2}^{2(\gamma+1)}
+\|\nabla v\|_{2}^{2(\gamma+1)}\right)\right],\label{2.3}
\end{align}
\begin{align}
 E(t):=E(u(t),v(t))=&\frac{1}{2}(\| u_{t}\| _{2}^{2}+\| v_{t}\|
_{2}^{2})+J(u(t),v(t)).
 \label{2.4}
\end{align}

As in \cite{2009a}, we can get
\begin{align}
E'(t)&=-(\|\nabla u_{t}\|_{2}^{2}+\|\nabla
v_{t}\|_{2}^{2})-\frac{1}{2}[g_{1}(t)\|\nabla
u\|_{2}^{2}+g_{2}(t)\|\nabla v\|_{2}^{2}]+\frac{1}{2}[(g_{1}'\circ
\nabla u)(t)+(g_{2}'\circ \nabla v)(t)] \nonumber \\
&\leq0, \quad
\forall\ t\geq 0. \label{2.5}
\end{align}
 Then we have
\begin{align}
E(t)+\int_{s}^{t}(\|\nabla u_{t}(\tau)\|_{2}^{2}+\|\nabla
v_{t}(\tau)\|_{2}^{2}){\rm d}\tau\leq E(s), \label{2.6}
\end{align}
for $0\leq s\leq t\leq T$.

Then, we give two lemmas which will be used throughout this work.

\begin{lemma}$($Sobolev-Poincar\'e inequality $\cite{18c}$$)$   \label{le2.1}
If $2\leq p\leq \frac{2N}{N-2},$ then
$$\|u\|_{p}\leq C_{p}\|\nabla u\|_{2,}$$
for $u\in H_{0}^{1}(\Omega)$ holds with some constants $C_{p}.$
\end{lemma}

\begin{lemma}$($See $\cite{17b}$$)$   \label{le2.2}
Let $h:[0,\infty)\rightarrow[0,\infty)$ be a non-increasing
function and assume that there exists a constant $r>0$ such that
$$\int_{t}^{\infty}h(s){\rm d}s\leq rh(t), \quad \forall\ t\in [0,\infty).$$
Then we have
$$h(t)\leq h(0)e^{1-\frac{t}{r}}, \quad \forall\ t\geq r.$$
\end{lemma}

 We now state a local existence theorem for system
\eqref{1.1}, whose proof follows the arguments in \cite{2006a,
2007a}:

\begin{theorem} \label{Th2.3}
Suppose that \eqref{2.1}, $($A1$)$ and $($A2$)$ hold, and that
$u_{0}, v_{0}\in H_{0}^{1}(\Omega)\cap H^{2}(\Omega),$ $u_{1},
v_{1}\in L^{2}(\Omega).$ Then problem \eqref{1.1} has a unique local
solution
$$u, v\in C([0,T]; H_{0}^{1}(\Omega)\cap
H^{2}(\Omega)),\quad u_{t}, v_{t}\in C([0,T]; L^{2}(\Omega))\cap
L^{2}([0,T]; H_{0}^{1}(\Omega)),$$
 for some $T>0.$
Moreover, at least one of the following statements is valid:
\begin{align*}
(1)\ &T=\infty,\\
(2)\  &\displaystyle{\lim_{t\rightarrow
T^{-}}}(\|u_{t}\|_{2}^{2}+\|v_{t}\|_{2}^{2}+\|\Delta
u\|_{2}^{2}+\|\Delta v\|_{2}^{2})=\infty.
\end{align*}
\end{theorem}

\section{Global existence and energy decay}

In this section, we consider the global existence and energy decay
of solutions for problem \eqref{1.1}. We first introduce two lemmas
which are essential in our proof.

\begin{lemma}{\rm(\cite[Lemma 3.2]{2010b})} \label{le3.1}
Assume that \eqref{2.1} holds. Then there exists $\eta>0$ such that for
any $(u,v)\in(H_{0}^{1}(\Omega)\cap
H^{2}(\Omega))\times(H_{0}^{1}(\Omega)\cap H^{2}(\Omega)),$ we have
\begin{align}
\|u+v\|_{2(p+2)}^{2(p+2)}+2\|uv\|_{p+2}^{p+2}\leq \eta(k_{1}\|\nabla
u\|_{2}^{2}+k_{2}\|\nabla v\|_{2}^{2})^{p+2}. \label{3.1}
\end{align}
\end{lemma}

In order to prove our result and for the sake of simplicity, we take
$a=b=1$ and introduce
\begin{align}
B=\eta^{\frac{1}{2(p+2)}}, \quad \alpha_{*}=B^{-\frac{p+2}{p+1}},
\quad E_{1}=\left(\frac{1}{2}-\frac{1}{2(p+2)}\right)\alpha_{*}^{2},
\label{3.2}
\end{align}
where $\eta$ is the optimal constant in \eqref{3.1}. Next, we will
state and prove a lemma which similar to the one introduced firstly
by Vitillaro in \cite{16b} to study a class of a single wave equation.

\begin{lemma}  \label{le3.2}
 Suppose that \eqref{2.1}, $($A1$)$ and $($A2$)$ hold. Let $(u,v)$ be the
 solution of system \eqref{1.1}. Assume further that $E(0)<E_{1}$
 and
\begin{align}
 (k_{1}\|\nabla
u_{0}\|_{2}^{2}+k_{2}\|\nabla v_{0}\|_{2}^{2})^{1/2}<\alpha_{*}.
\label{3.3}
\end{align}
Then
\begin{align}
\left(k_{1}\|\nabla u(t)\|_{2}^{2}+k_{2}\|\nabla
v(t)\|_{2}^{2}+(g_{1}\circ \nabla u)(t)+(g_{2}\circ \nabla
v)(t)\right)^{1/2}<\alpha_{*}, \label{3.4}
\end{align}
for all $t\in[0,T).$
\end{lemma}

{\bf Proof.}\ We first note that, by \eqref{2.4}, \eqref{3.1} and
the definition of $B,$ we have
\begin{align}
E(t)\geq&\frac{1}{2}\left[k_{1}\|\nabla u(t)\|_{2}^{2}+k_{2}\|\nabla
v(t)\|_{2}^{2}+(g_{1}\circ \nabla u)(t)+(g_{2}\circ \nabla
v)(t)\right] \nonumber\\
&-\frac{1}{2(p+2)}\left(\|u+v\|_{2(p+2)}^{2(p+2)}+2\|uv\|_{p+2}^{p+2}\right) \nonumber\\
\geq&\frac{1}{2}\left[k_{1}\|\nabla u(t)\|_{2}^{2}+k_{2}\|\nabla
v(t)\|_{2}^{2}+(g_{1}\circ \nabla u)(t)+(g_{2}\circ \nabla
v)(t)\right] \nonumber\\
&-\frac{B^{2(p+2)}}{2(p+2)}\left(k_{1}\|\nabla
u(t)\|_{2}^{2}+k_{2}\|\nabla v(t)\|_{2}^{2}\right)^{p+2} \nonumber\\
\geq&\frac{1}{2}\alpha^{2}-\frac{B^{2(p+2)}}{2(p+2)}\alpha^{2(p+2)}\triangleq
G(\alpha), \label{3.5}
\end{align}
where
$$\alpha=\left(k_{1}\|\nabla u(t)\|_{2}^{2}+k_{2}\|\nabla
v(t)\|_{2}^{2}+(g_{1}\circ \nabla u)(t)+(g_{2}\circ \nabla
v)(t)\right)^{1/2}.$$ It is easy to verify that $G(\alpha)$ is
increasing in $(0,\alpha_{*})$, decreasing in $(\alpha_{*},\infty),$
and that $G(\alpha)\rightarrow-\infty,$ as
$\alpha\rightarrow\infty,$ and
$$G(\alpha)_{\max}=G(\alpha_{*})=\frac{1}{2}\alpha_{*}^{2}-\frac{B^{2(p+2)}}{2(p+2)}\alpha_{*}^{2(p+2)}=E_{1},$$
where $\alpha_{*}$ is given in \eqref{3.2}.

 Now we establish
\eqref{3.4} by contradiction. First we assume that \eqref{3.4} is
not true over $[0,T)$, then it follows from the continuity of
$(u(t), v(t))$ that there exists $t_{0}\in (0,T)$ such that
$$\left(k_{1}\|\nabla u(t_{0})\|_{2}^{2}+k_{2}\|\nabla
v(t_{0})\|_{2}^{2}+(g_{1}\circ \nabla u)(t_{0})+(g_{2}\circ \nabla
v)(t_{0})\right)^{1/2}=\alpha_{*}.$$ By \eqref{3.5}, we see that
\begin{align*}
E(t_{0})&\geq G\left[\left(k_{1}\|\nabla
u(t_{0})\|_{2}^{2}+k_{2}\|\nabla v(t_{0})\|_{2}^{2}+(g_{1}\circ
\nabla u)(t_{0})+(g_{2}\circ \nabla v)(t_{0})\right)^{1/2}\right] \\
&=G(\alpha_{*})=E_{1}.
\end{align*}
This is impossible since $E(t)\leq E(0)<E_{1},$  $t\geq0.$ Hence
\eqref{3.4} is established. \quad $\Box$

Now we are in a position to state and prove our main result.
\begin{theorem} $($\emph{Global existence and energy decay}$)$ \label{Th3.3}
Assume that \eqref{2.1}, $($A1$)$ and $($A2$)$ hold. If the initial data
$(u_{0}, u_{1})\in(H_{0}^{1}(\Omega)\cap H^{2}(\Omega))\times
L^{2}(\Omega),$ $(v_{0}, v_{1})\in(H_{0}^{1}(\Omega)\cap
H^{2}(\Omega))\times L^{2}(\Omega),$ satisfy $E(0)<E_{1}$ and
\begin{align}
(k_{1}\|\nabla u_{0}\|_{2}^{2}+k_{2}\|\nabla
v_{0}\|_{2}^{2})^{1/2}<\alpha_{*}, \label{3.6}
\end{align}
where the constants $\alpha_{*}$, $E_{1}$ are defined in
\eqref{3.2}, then the corresponding solution to system \eqref{1.1}
globally exists, i.e., $T=\infty.$

Moreover, if $m_{0}-k>0$ is sufficiently small such that
\begin{align}
1-\eta\left(\frac{2(p+2)}{p+1}E(0)\right)^{p+1}-\frac{5(m_{0}-k)(p+2)}{2k(p+1)}>0,
\label{3.7}
\end{align}
 then we have the following decay estimates
$$E(t)\leq E(0)e^{1-\lambda C_{0}^{-1}t}$$
for every $t\geq \lambda C_{0}^{-1},$ where $C_{0}$ is some positive
constant.
\end{theorem}

{\bf Proof.}\ First, we prove that $T=\infty.$ Since $E(0)<E_{1}$
and
 $$(k_{1}\|\nabla u_{0}\|_{2}^{2}+k_{2}\|\nabla
v_{0}\|_{2}^{2})^{1/2}<\alpha_{*},$$ it follows from Lemma
\ref{le3.2} that
\begin{align*}
&k_{1}\|\nabla u(t)\|_{2}^{2}+k_{2}\|\nabla v(t)\|_{2}^{2} \\
\leq& k_{1}\|\nabla u(t)\|_{2}^{2}+k_{2}\|\nabla
v(t)\|_{2}^{2}+(g_{1}\circ \nabla u)(t)+(g_{2}\circ \nabla v)(t) \\
<&\alpha_{*}^{2}=\eta^{-\frac{1}{p+1}},
\end{align*}
which implies that
\begin{align*}
I(t)\geq&\left(m_{0}-\int_{0}^{\infty}g_{1}(s){\rm
d}s\right)\|\nabla
u\|_{2}^{2}+\left(m_{0}-\int_{0}^{\infty}g_{2}(s){\rm
d}s\right)\|\nabla v\|_{2}^{2}
\\ &+(g_{1}\circ \nabla u)(t)+(g_{2}\circ \nabla
v)(t)-2(p+2)\int_{\Omega}F(u,v){\rm d}x \\
 \geq& k_{1}\|\nabla u(t)\|_{2}^{2}+k_{2}\|\nabla
v(t)\|_{2}^{2}-2(p+2)\int_{\Omega}F(u,v){\rm d}x \\
=&k_{1}\|\nabla u(t)\|_{2}^{2}+k_{2}\|\nabla
v(t)\|_{2}^{2}-\left(\|u+v\|_{2(p+2)}^{2(p+2)}+2\|uv\|_{p+2}^{p+2}\right)
\\
\geq& k_{1}\|\nabla u(t)\|_{2}^{2}+k_{2}\|\nabla
v(t)\|_{2}^{2}-\eta\left(k_{1}\|\nabla u(t)\|_{2}^{2}+k_{2}\|\nabla
v(t)\|_{2}^{2}\right)^{p+2} \\
=&(k_{1}\|\nabla u(t)\|_{2}^{2}+k_{2}\|\nabla
v(t)\|_{2}^{2})\left[1-\eta\left(k_{1}\|\nabla
u(t)\|_{2}^{2}+k_{2}\|\nabla
v(t)\|_{2}^{2}\right)^{p+1}\right]\geq0,
\end{align*}
\nopagebreak[1] for $t\in[0,T),$ where we have used \eqref{3.1}.
Further, by \eqref{2.2} and \eqref{2.3}, we have
\begin{align}
J(t)\geq&\frac{1}{2}\left[\left(m_{0}-\int_{0}^{t} g_{1}(s){\rm
d}s\right)\|\nabla u\|_{2}^{2}+\left(m_{0}-\int_{0}^{t} g_{2}(s){\rm
d}s\right)\|\nabla v\|_{2}^{2} +(g_{1}\circ \nabla u)(t)+(g_{2}\circ
\nabla
v)(t)\right] \nonumber\\
&-\int_{\Omega}F(u,v){\rm d}x-\frac{1}{2(p+2)}I(t)+\frac{1}{2(p+2)}I(t) \nonumber\\
\geq&\left(\frac{1}{2}-\frac{1}{2(p+2)}\right)\left[\left(m_{0}-\int_{0}^{t}
g_{1}(s){\rm d}s\right)\|\nabla u\|_{2}^{2}+\left(m_{0}-\int_{0}^{t}
g_{2}(s){\rm d}s\right)\|\nabla v\|_{2}^{2}\right]\nonumber \\
&+\left(\frac{1}{2}-\frac{1}{2(p+2)}\right)\left[(g_{1}\circ \nabla
u)(t)+(g_{2}\circ \nabla v)(t)\right]+\frac{1}{2(p+2)}I(t) \nonumber\\
\geq&\frac{p+1}{2(p+2)}\left[k_{1}\|\nabla
u(t)\|_{2}^{2}+k_{2}\|\nabla v(t)\|_{2}^{2}+(g_{1}\circ \nabla
u)(t)+(g_{2}\circ \nabla v)(t)\right]+\frac{1}{2(p+2)}I(t)
\nonumber\\ \geq&0, \label{3.8}
\end{align}
from \eqref{3.8} and \eqref{2.4}, we deduce that
\begin{align}
k\left(\|\nabla u(t)\|_{2}^{2}+\|\nabla v(t)\|_{2}^{2}\right)\leq
k_{1}\|\nabla u(t)\|_{2}^{2}+k_{2}\|\nabla
v(t)\|_{2}^{2}\leq\frac{2(p+2)}{p+1}J(t)\leq&\frac{2(p+2)}{p+1}E(t)
. \label{3.9}
\end{align}

Multiplying the first equation of system \eqref{1.1} by $-2\Delta
u,$ and integrating it over $\Omega,$ we get
\begin{align*}
&\frac{d}{dt}\left\{\|\Delta u\|_{2}^{2}-2\int_{\Omega}u_{t}\Delta
u{\rm d}x\right\}+2M(\|\nabla u\|_{2}^{2})\|\Delta u\|_{2}^{2} \\
\leq&2\|\nabla u_{t}\|_{2}^{2}-2\int_{\Omega}f_{1}(u,v)\Delta u{\rm
d}x+2\int_{0}^{t}g_{1}(t-s)\int_{\Omega}\Delta u(s)\Delta u(t){\rm
d}x{\rm d}s
\end{align*}
By Young's inequality, we have
\begin{align*}
2\int_{0}^{t}g_{1}(t-s)\int_{\Omega}\Delta u(s)\Delta u(t){\rm
d}x{\rm d}s\leq2\theta\|\Delta
u(t)\|_{2}^{2}+\frac{\|g_{1}\|_{L^{1}}}{2\theta}\int_{0}^{t}g_{1}(t-s)\|\Delta
u(s)\|_{2}^{2}{\rm d}s.
\end{align*}
Thus, we obtain
\begin{align}
&\frac{d}{dt}\left\{\|\Delta u\|_{2}^{2}-2\int_{\Omega}u_{t}\Delta
u{\rm d}x\right\}+\left[2M(\|\nabla u\|_{2}^{2})-2\theta\right]\|\Delta u\|_{2}^{2} \nonumber\\
\leq&2\|\nabla
u_{t}\|_{2}^{2}+\frac{\|g_{1}\|_{L^{1}}}{2\theta}\int_{0}^{t}g_{1}(t-s)\|\Delta
u(s)\|_{2}^{2}{\rm d}s-2\int_{\Omega}f_{1}(u,v)\Delta u{\rm d}x,
\label{3.10}
\end{align}
where $0<\theta\leq\frac{\|g_{0}\|_{L^{1}}}{2}.$ Similarly, we also
have
\begin{align}
&\frac{d}{dt}\left\{\|\Delta v\|_{2}^{2}-2\int_{\Omega}v_{t}\Delta
v{\rm d}x\right\}+\left[2M(\|\nabla v\|_{2}^{2})-2\theta\right]\|\Delta v\|_{2}^{2} \nonumber\\
\leq&2\|\nabla
v_{t}\|_{2}^{2}+\frac{\|g_{2}\|_{L^{1}}}{2\theta}\int_{0}^{t}g_{2}(t-s)\|\Delta
v(s)\|_{2}^{2}{\rm d}s-2\int_{\Omega}f_{2}(u,v)\Delta v{\rm d}x.
\label{3.11}
\end{align}
Now, combining \eqref{3.10} with \eqref{3.11}, we get
\begin{align}
&\frac{d}{dt}\left\{\|\Delta u\|_{2}^{2}-2\int_{\Omega}u_{t}\Delta
u{\rm d}x+\|\Delta v\|_{2}^{2}-2\int_{\Omega}v_{t}\Delta v{\rm d}x\right\}  \nonumber\\
&+2\left[M(\|\nabla u\|_{2}^{2})-\theta\right]\|\Delta u\|_{2}^{2}+
2\left[M(\|\nabla v\|_{2}^{2})-\theta\right]\|\Delta v\|_{2}^{2} \nonumber\\
\leq&2(\|\nabla u_{t}\|_{2}^{2}+\|\nabla
v_{t}\|_{2}^{2})+\frac{\|g_{1}\|_{L^{1}}}{2\theta}\int_{0}^{t}g_{1}(t-s)\|\Delta
u(s)\|_{2}^{2}{\rm d}s \nonumber\\
&+\frac{\|g_{2}\|_{L^{1}}}{2\theta}\int_{0}^{t}g_{2}(t-s)\|\Delta
v(s)\|_{2}^{2}{\rm d}s-2\int_{\Omega}(f_{1}(u,v)\Delta
u+f_{2}(u,v)\Delta v){\rm d}x. \label{3.12}
\end{align}
Multiplying \eqref{3.12} by $\varepsilon,$ $0<\varepsilon\leq1,$ and
multiplying \eqref{2.5} by $\delta,$ $\delta$ is some positive
constant, and then summing up, we deduce
\begin{align}
&\frac{d}{dt}E^{*}(t)+(\delta-2\varepsilon)(\|\nabla
u_{t}\|_{2}^{2}+\|\nabla
v_{t}\|_{2}^{2})+2\varepsilon\left[M(\|\nabla
u\|_{2}^{2})-\theta\right]\|\Delta
u\|_{2}^{2}+2\varepsilon\left[M(\|\nabla
v\|_{2}^{2})-\theta\right]\|\Delta v\|_{2}^{2}  \nonumber\\
\leq&\varepsilon\frac{\|g_{1}\|_{L^{1}}}{2\theta}\int_{0}^{t}g_{1}(t-s)\|\Delta
u(s)\|_{2}^{2}{\rm
d}s+\varepsilon\frac{\|g_{2}\|_{L^{1}}}{2\theta}\int_{0}^{t}g_{2}(t-s)\|\Delta
v(s)\|_{2}^{2}{\rm d}s
\nonumber\\&-2\varepsilon\int_{\Omega}(f_{1}(u,v)\Delta
u+f_{2}(u,v)\Delta v){\rm d}x, \label{3.13}
\end{align}
where
$$E^{*}(t)=\delta E(t)+\varepsilon(\|\Delta u\|_{2}^{2}+\|\Delta v\|_{2}^{2})-2\varepsilon\left(\int_{\Omega}u_{t}\Delta
u{\rm d}x+\int_{\Omega}v_{t}\Delta v{\rm d}x\right).$$ By Young's
inequality, we get $$\left|2\varepsilon\int_{\Omega}u_{t}\Delta
u{\rm
d}x\right|\leq2\varepsilon\|u_{t}\|_{2}^{2}+\frac{\varepsilon}{2}\|\Delta
u\|_{2}^{2}$$ and $$\left|2\varepsilon\int_{\Omega}v_{t}\Delta v{\rm
d}x\right|\leq2\varepsilon\|v_{t}\|_{2}^{2}+\frac{\varepsilon}{2}\|\Delta
v\|_{2}^{2}.$$ Noting that $J(t)\geq0$ by \eqref{3.8}, then, by
\eqref{2.4}, we have
\begin{align}
\|u_{t}\|_{2}^{2}+\|v_{t}\|_{2}^{2}\leq2E(t). \label{3.14}
\end{align}
Therefore, choosing $\delta=5\varepsilon,$ we observe that
\begin{align}
E^{*}(t)\geq\frac{\varepsilon}{2}\left(\|u_{t}\|_{2}^{2}+\|v_{t}\|_{2}^{2}+\|\Delta
u\|_{2}^{2}+\|\Delta v\|_{2}^{2}\right). \label{3.15}
\end{align}
Moreover, by H$\ddot{\text{o}}$lder's inequality, Lemma \ref{le2.1}
and \eqref{3.9}, we see that
\begin{align}
&\left|\int_{\Omega}f_{1}(u,v)\Delta u{\rm d}x\right| \nonumber\\
\leq&\int_{\Omega}\left(|u+v|^{2p+3}+|u|^{p+1}|v|^{p+2}\right)\Delta
u{\rm d}x  \nonumber\\
\leq&C\int_{\Omega}\left(|u|^{2p+3}+|v|^{2p+3}+|u|^{p+1}|v|^{p+2}\right)\Delta
u{\rm d}x  \nonumber\\
\leq&C\left[\int_{\Omega}\left(|u|^{2p+3}+|v|^{2p+3}+|u|^{p+1}|v|^{p+2}\right)^{2}{\rm
d}x\right]^{\frac{1}{2}}\|\Delta
u\|_{2} \nonumber\\
\leq&C\left[\int_{\Omega}\left(|u|^{2(2p+3)}+|v|^{2(2p+3)}+|u|^{2(p+1)}|v|^{2(p+2)}\right){\rm
d}x\right]^{\frac{1}{2}}\|\Delta
u\|_{2} \nonumber\\
\leq&C\left[\|u\|_{2(2p+3)}^{2p+3}+\|v\|_{2(2p+3)}^{2p+3}+\|u\|_{2(2p+3)}^{p+1}\|v\|_{2(2p+3)}^{p+2}\right]\|\Delta
u\|_{2} \nonumber\\
\leq&C\left[\|\nabla u\|_{2}^{2p+3}+\|\nabla v\|_{2}^{2p+3}+\|\nabla
u\|_{2}^{p+1}\|\nabla v\|_{2}^{p+2}\right]\|\Delta
u\|_{2} \nonumber\\
\leq&C_{*}\|\Delta u\|_{2}, \label{3.16}
\end{align}
here we denote by $C>0$ a generic constant that may vary even from
line to line within the same formula and
$C_{*}=3C\left(\frac{2(p+2)}{k(p+1)}E(0)\right)^{(2p+3)/2}.$
Similarly, we have
\begin{align}
\left|\int_{\Omega}f_{2}(u,v)\Delta v{\rm d}x\right|\leq
C_{*}\|\Delta v\|_{2}. \label{3.17}
\end{align}
Combining \eqref{3.15}--\eqref{3.17}, we deduce
\begin{align}
2\varepsilon\left|\int_{\Omega}f_{1}(u,v)\Delta u{\rm
d}x+\int_{\Omega}f_{2}(u,v)\Delta v{\rm d}x\right|\leq2\varepsilon
C_{*}(\|\Delta u\|_{2}+\|\Delta v\|_{2})\leq
4\sqrt{2\varepsilon}C_{*}E^{*}(t)^{\frac{1}{2}}. \label{3.18}
\end{align}
Substituting \eqref{3.18} into \eqref{3.13}, and then integrating it
over $(0,t),$ we obtain
\begin{align}
&E^{*}(t)+2\varepsilon\left(m_{0}-\theta-\frac{\|g_{0}\|_{L^{1}}^{2}}{4\theta}\right)\int_{0}^{t}(\|\Delta
u(s)\|_{2}^{2}+\|\Delta v(s)\|_{2}^{2}){\rm d}s \nonumber\\
 \leq&
E^{*}(0)+4\sqrt{2\varepsilon}C_{*}\int_{0}^{t}E^{*}(s)^{\frac{1}{2}}{\rm
d}s. \label{3.19}
\end{align}
Taking $\theta=\frac{\|g_{0}\|_{L^{1}}}{2}$ in \eqref{3.19}, and
then by Gronwall's Lemma, we deduce
$$E^{*}(t)\leq\left(2\sqrt{2\varepsilon}C_{*}t+E^{*}(0)^{\frac{1}{2}}\right)^{2},$$
for any $t\geq0.$ Therefore, by \eqref{3.15} and Theorem
\ref{Th2.3}, we have $T=\infty$ as long as $\varepsilon$ is fixed.

Next, we want to derive the decay rate of energy function for system
\eqref{1.1}. Multiplying the first equation of system \eqref{1.1} by
$u$ and the second equation of system \eqref{1.1} by $v,$
integrating them over $\Omega\times[t_{1}, t_{2}] (0\leq t_{1}\leq
t_{2}),$ using integration by parts and summing up, we have
\begin{align*}
&\int_{\Omega}u_{t}u{\rm
d}x\bigg|_{t_{1}}^{t_{2}}-\int_{t_{1}}^{t_{2}}\|u_{t}\|_{2}^{2}{\rm
d}t +\int_{\Omega}v_{t}v{\rm
d}x\bigg|_{t_{1}}^{t_{2}}-\int_{t_{1}}^{t_{2}}\|v_{t}\|_{2}^{2}{\rm
d}t +\int_{t_{1}}^{t_{2}}M(\|\nabla u\|_{2}^{2})\|\nabla
u\|_{2}^{2}{\rm d}t \\
&+\int_{t_{1}}^{t_{2}}M(\|\nabla v\|_{2}^{2})\|\nabla
v\|_{2}^{2}{\rm
d}t+\int_{t_{1}}^{t_{2}}\int_{\Omega}\int_{0}^{t}g_{1}(t-s)\Delta
u(s)u(t){\rm d}s{\rm d}x{\rm d}t \\
&+\int_{t_{1}}^{t_{2}}\int_{\Omega}\int_{0}^{t}g_{2}(t-s)\Delta
v(s)v(t){\rm d}s{\rm d}x{\rm
d}t+\int_{t_{1}}^{t_{2}}\int_{\Omega}\nabla u_{t}\nabla u{\rm
d}x{\rm d}t+\int_{t_{1}}^{t_{2}}\int_{\Omega}\nabla v_{t}\nabla
v{\rm d}x{\rm d}t \\
=&\int_{t_{1}}^{t_{2}}\int_{\Omega}[f_{1}(u,v) u+f_{2}(u,v) v]{\rm
d}x{\rm d}t.
\end{align*}
Consequently,
\begin{align*}
&\int_{t_{1}}^{t_{2}}M(\|\nabla u\|_{2}^{2})\|\nabla u\|_{2}^{2}{\rm
d}t+\int_{t_{1}}^{t_{2}}M(\|\nabla v\|_{2}^{2})\|\nabla
v\|_{2}^{2}{\rm d}t-2(p+2)\int_{t_{1}}^{t_{2}}\int_{\Omega}F(u,v){\rm d}x{\rm d}t \\
=&-\int_{\Omega}u_{t}u{\rm
d}x\bigg|_{t_{1}}^{t_{2}}-\int_{\Omega}v_{t}v{\rm
d}x\bigg|_{t_{1}}^{t_{2}} -\int_{t_{1}}^{t_{2}}\int_{\Omega}\nabla
u_{t}\nabla u{\rm d}x{\rm
d}t-\int_{t_{1}}^{t_{2}}\int_{\Omega}\nabla v_{t}\nabla
v{\rm d}x{\rm d}t+\int_{t_{1}}^{t_{2}}\|u_{t}\|_{2}^{2}{\rm d}t \\
&+\int_{t_{1}}^{t_{2}}\|v_{t}\|_{2}^{2}{\rm d}t
-\int_{t_{1}}^{t_{2}}\int_{\Omega}\int_{0}^{t}g_{1}(t-s)\Delta
u(s)u(t){\rm d}s{\rm d}x{\rm d}t
-\int_{t_{1}}^{t_{2}}\int_{\Omega}\int_{0}^{t}g_{2}(t-s)\Delta
v(s)v(t){\rm d}s{\rm d}x{\rm d}t.
\end{align*}
It follows from \eqref{2.4} that
\begin{align*}
&2\int_{t_{1}}^{t_{2}}E(t){\rm
d}t+\int_{t_{1}}^{t_{2}}\left[M(\|\nabla u\|_{2}^{2})\|\nabla
u\|_{2}^{2}+M(\|\nabla v\|_{2}^{2})\|\nabla
v\|_{2}^{2}\right]{\rm d}t-2(p+2)\int_{t_{1}}^{t_{2}}\int_{\Omega}F(u,v){\rm d}x{\rm d}t \\
=&2\int_{t_{1}}^{t_{2}}(\|u_{t}\|_{2}^{2}+\|v_{t}\|_{2}^{2}){\rm
d}t+\int_{t_{1}}^{t_{2}}\left[\left(m_{0}-\int_{0}^{t} g_{1}(s){\rm
d}s\right)\|\nabla u\|_{2}^{2}+\left(m_{0}-\int_{0}^{t}
g_{2}(s){\rm d}s\right)\|\nabla v\|_{2}^{2}\right]{\rm d}t \\
&-2\int_{t_{1}}^{t_{2}}\int_{\Omega}F(u,v){\rm d}x{\rm
d}t-\int_{\Omega}u_{t}u{\rm d}x\bigg|_{t_{1}}^{t_{2}}-
\int_{\Omega}v_{t}v{\rm d}x\bigg|_{t_{1}}^{t_{2}}
-\int_{t_{1}}^{t_{2}}\int_{\Omega}(\nabla u_{t}\nabla u+\nabla v_{t}\nabla v){\rm d}x{\rm d}t \\
& -\int_{t_{1}}^{t_{2}}\int_{\Omega}\int_{0}^{t}g_{1}(t-s)\Delta
u(s)u(t){\rm d}s{\rm d}x{\rm d}t
-\int_{t_{1}}^{t_{2}}\int_{\Omega}\int_{0}^{t}g_{2}(t-s)\Delta
v(s)v(t){\rm d}s{\rm d}x{\rm d}t \\
&+\int_{t_{1}}^{t_{2}}\left[(g_{1}\circ \nabla u)(t)+(g_{2}\circ
\nabla v)(t)+\frac{m_{1}}{\gamma+1}\left(\|\nabla
u\|_{2}^{2(\gamma+1)}+\|\nabla
v\|_{2}^{2(\gamma+1)}\right)\right]{\rm d}t,
\end{align*}
\nopagebreak [2] and then we arrive at
\begin{align}
&2\int_{t_{1}}^{t_{2}}E(t){\rm
d}t-2(p+1)\int_{t_{1}}^{t_{2}}\int_{\Omega}F(u,v){\rm d}x{\rm d}t
\nonumber \\
 \leq&-\int_{t_{1}}^{t_{2}}\left[\int_{0}^{t}g_{1}(s){\rm d}s\|\nabla
u\|_{2}^{2}+\int_{0}^{t}g_{2}(s){\rm d}s\|\nabla
v\|_{2}^{2}\right]{\rm d}t+\int_{t_{1}}^{t_{2}}\left[(g_{1}\circ
\nabla
u)(t)+(g_{2}\circ \nabla v)(t)\right]{\rm d}t \nonumber\\
&-\int_{\Omega}(u_{t}u+v_{t}v){\rm d}x\bigg|_{t_{1}}^{t_{2}}+
2\int_{t_{1}}^{t_{2}}(\|u_{t}\|_{2}^{2}+\|v_{t}\|_{2}^{2}){\rm
d}t-\int_{t_{1}}^{t_{2}}\int_{\Omega}(\nabla
u_{t}\nabla u+\nabla v_{t}\nabla v){\rm d}x{\rm d}t  \nonumber\\
&-\int_{t_{1}}^{t_{2}}\int_{\Omega}\int_{0}^{t}g_{1}(t-s)\Delta
u(s)u(t){\rm d}s{\rm d}x{\rm d}t
-\int_{t_{1}}^{t_{2}}\int_{\Omega}\int_{0}^{t}g_{2}(t-s)\Delta
v(s)v(t){\rm d}s{\rm d}x{\rm d}t.  \label{3.20}
\end{align}
For the sixth term on the right-hand side of \eqref{3.20}, we have
\begin{align}
&-\int_{\Omega}\int_{0}^{t}g_{1}(t-s)\Delta u(s)u(t){\rm d}s{\rm
d}x{\rm d}t=\int_{\Omega}\int_{0}^{t}g_{1}(t-s)\nabla u(s)\nabla
u(t){\rm d}s{\rm d}x \nonumber\\
=&\frac{1}{2}\left[\int_{0}^{t}g_{1}(t-s)\|\nabla u(t)\|_{2}^{2}{\rm
d}s+\int_{0}^{t}g_{1}(t-s)\|\nabla u(s)\|_{2}^{2}{\rm
d}s-(g_{1}\circ \nabla u)(t)\right]. \label{3.21}
\end{align}
 Similarly,
\begin{align}
&-\int_{\Omega}\int_{0}^{t}g_{2}(t-s)\Delta v(s)v(t){\rm d}s{\rm
d}x{\rm d}t
\nonumber\\
=&\frac{1}{2}\left[\int_{0}^{t}g_{2}(t-s)\|\nabla v(t)\|_{2}^{2}{\rm
d}s+\int_{0}^{t}g_{2}(t-s)\|\nabla v(s)\|_{2}^{2}{\rm
d}s-(g_{2}\circ \nabla v)(t)\right]. \label{3.22}
\end{align}
By \eqref{3.20}, \eqref{3.21} and \eqref{3.22}, we get
\begin{align}
&2\int_{t_{1}}^{t_{2}}E(t){\rm d}t-2(p+1)\int_{t_{1}}^{t_{2}}\int_{\Omega}F(u,v){\rm d}x{\rm d}t \nonumber\\
\leq&-\frac{1}{2}\int_{t_{1}}^{t_{2}}\left[\int_{0}^{t}g_{1}(s){\rm
d}s\|\nabla u(t)\|_{2}^{2}+\int_{0}^{t}g_{2}(s){\rm d}s\|\nabla
v(t)\|_{2}^{2}\right]{\rm
d}t+\frac{1}{2}\int_{t_{1}}^{t_{2}}\left[(g_{1}\circ \nabla
u)(t)+(g_{2}\circ \nabla v)(t)\right]{\rm d}t \nonumber\\
&-\int_{\Omega}(u_{t}u+v_{t}v){\rm d}x\bigg|_{t_{1}}^{t_{2}}+
2\int_{t_{1}}^{t_{2}}(\|u_{t}\|_{2}^{2}+\|v_{t}\|_{2}^{2}){\rm
d}t-\int_{t_{1}}^{t_{2}}\int_{\Omega}(\nabla
u_{t}\nabla u+\nabla v_{t}\nabla v){\rm d}x{\rm d}t  \nonumber\\
&+\frac{1}{2}\int_{t_{1}}^{t_{2}}\int_{0}^{t}g_{1}(t-s)\|\nabla
u(s)\|_{2}^{2}{\rm d}s{\rm
d}t+\frac{1}{2}\int_{t_{1}}^{t_{2}}\int_{0}^{t}g_{2}(t-s)\|\nabla
v(s)\|_{2}^{2}{\rm d}s{\rm d}t \nonumber\\
\leq&-\int_{\Omega}(u_{t}u+v_{t}v){\rm d}x\bigg|_{t_{1}}^{t_{2}}+
2\int_{t_{1}}^{t_{2}}(\|u_{t}\|_{2}^{2}+\|v_{t}\|_{2}^{2}){\rm
d}t-\int_{t_{1}}^{t_{2}}\int_{\Omega}(\nabla
u_{t}\nabla u+\nabla v_{t}\nabla v){\rm d}x{\rm d}t  \nonumber\\
&+\frac{1}{2}\int_{t_{1}}^{t_{2}}\left[\int_{0}^{t}g_{1}(t-s)\|\nabla
u(s)\|_{2}^{2}{\rm d}s+\int_{0}^{t}g_{2}(t-s)\|\nabla
v(s)\|_{2}^{2}{\rm d}s\right]{\rm d}t \nonumber\\
&+\frac{1}{2}\int_{t_{1}}^{t_{2}}\left[(g_{1}\circ \nabla
u)(t)+(g_{2}\circ \nabla v)(t)\right]{\rm d}t \nonumber\\
=&I_{1}+I_{2}-I_{3}+I_{4}+I_{5}. \label{3.23}
\end{align}

In what follows we will estimate $I_{1},\ldots, I_{5}$ in
\eqref{3.23}. Firstly, by H$\ddot{\text{o}}$lder's inequality,
Young's inequality and Lemma \ref{le2.1}, we have
\begin{align*}
& \int_{\Omega}|u(t)u_{t}(t)|{\rm d}x+\int_{\Omega}|v(t)v_{t}(t)|{\rm d}x  \\
\leq&\frac{1}{2}\|u(t)\|_{2}^{2}+\frac{1}{2}\|u_{t}(t)\|_{2}^{2}+\frac{1}{2}\|v(t)\|_{2}^{2}
+\frac{1}{2}\|v_{t}(t)\|_{2}^{2}  \\
\leq&\frac{C_{p}^{2}}{2}\|\nabla
u(t)\|_{2}^{2}+\frac{C_{p}^{2}}{2}\|\nabla
v(t)\|_{2}^{2}+\frac{1}{2}\|u_{t}(t)\|_{2}^{2}+\frac{1}{2}\|v_{t}(t)\|_{2}^{2}.
\end{align*}
It also follows from \eqref{3.9}, \eqref{2.4} and \eqref{2.5} that
$k(\|\nabla u(t)\|_{2}^{2}+\|\nabla
v(t)\|_{2}^{2})\leq\frac{2(p+2)}{p+1}E(t),$
$\|u_{t}(t)\|_{2}^{2}+\|v_{t}(t)\|_{2}^{2}\leq2E(t)$
 and $E(t)$ is a non-increasing function. Therefore, we have
\begin{align}
I_{1}\leq\int_{\Omega}|u(t)u_{t}(t)|{\rm
d}x\bigg|_{t_{1}}^{t_{2}}+\int_{\Omega}|v(t)v_{t}(t)|{\rm
d}x\bigg|_{t_{1}}^{t_{2}} \leq2C_{1}E(t_{1}), \label{3.24}
\end{align}
where $C_{1}=\frac{p+2}{k(p+1)}C_{p}^{2}+1.$

For $I_{2}$ in \eqref{3.23}, applying $\|\nabla
u_{t}\|_{2}^{2}+\|\nabla v_{t}\|_{2}^{2}\leq-E'(t)$ from
\eqref{2.5}, we have
\begin{align}
I_{2}\leq2C_{p}^{2}\int_{t_{1}}^{t_{2}}(\|\nabla
u_{t}(t)\|_{2}^{2}+\|\nabla v_{t}(t)\|_{2}^{2}){\rm
d}t\leq2C_{p}^{2}E(t_{1}),  \label{3.25}
\end{align}

We also have the following estimate
\begin{align}
I_{3}=&\int_{t_{1}}^{t_{2}}\int_{\Omega}\nabla u(t)\nabla
u_{t}(t){\rm d}x{\rm d}t+\int_{t_{1}}^{t_{2}}\int_{\Omega}\nabla
v(t)\nabla
v_{t}(t){\rm d}x{\rm d}t \nonumber\\
=&\frac{1}{2}\int_{t_{1}}^{t_{2}}\frac{d}{{\rm d}t}\|\nabla
u(t)\|_{2}^{2}{\rm d}t+\frac{1}{2}\int_{t_{1}}^{t_{2}}\frac{d}{{\rm
d}t}\|\nabla
v(t)\|_{2}^{2}{\rm d}t \nonumber\\
=&\frac{1}{2}\left(\|\nabla u(t_{2})\|_{2}^{2}-\|\nabla
u(t_{1})\|_{2}^{2}\right)+\frac{1}{2}\left(\|\nabla
v(t_{2})\|_{2}^{2}-\|\nabla
v(t_{1})\|_{2}^{2}\right) \nonumber\\
\leq&\frac{2(p+2)}{k(p+1)}E(t_{1})\triangleq C_{3}E(t_{1}).
\label{3.26}
\end{align}
To estimate $I_{4}$, using Young's inequality for convolution
$\|\phi\ast\psi\|_{q}\leq\|\phi\|_{r}\|\psi\|_{s},$ with
$\frac{1}{q}=\frac{1}{r}+\frac{1}{s}-1,$ $1\leq q, r, s\leq\infty,$
noting that if $q=1,$ then $r=1$ and $s=1,$ we get
\begin{align}
\int_{t_{1}}^{t_{2}}\int_{0}^{t}g_{1}(t-s)\|\nabla
u(s)\|_{2}^{2}{\rm d}s{\rm d}t\leq\int_{t_{1}}^{t_{2}}g_{1}(t){\rm
d}t\int_{t_{1}}^{t_{2}}\|\nabla u(t)\|_{2}^{2}{\rm d}t
 \leq(m_{0}-k_{1})\int_{t_{1}}^{t_{2}}\|\nabla
u(t)\|_{2}^{2}{\rm d}t, \label{3.27}
\end{align}
and
\begin{align}
\int_{t_{1}}^{t_{2}}\int_{0}^{t}g_{2}(t-s)\|\nabla
v(s)\|_{2}^{2}{\rm d}s{\rm d}t\leq\int_{t_{1}}^{t_{2}}g_{2}(t){\rm
d}t\int_{t_{1}}^{t_{2}}\|\nabla v(t)\|_{2}^{2}{\rm d}t
 \leq(m_{0}-k_{2})\int_{t_{1}}^{t_{2}}\|\nabla
v(t)\|_{2}^{2}{\rm d}t. \label{3.28}
\end{align}
Hence, by \eqref{3.9}, \eqref{3.27} and \eqref{3.28}, we obtain
\begin{align}
2I_{4}=&\int_{t_{1}}^{t_{2}}\int_{0}^{t}g_{1}(t-s)\|\nabla
u(s)\|_{2}^{2}{\rm d}s{\rm
d}t+\int_{t_{1}}^{t_{2}}\int_{0}^{t}g_{2}(t-s)\|\nabla
v(s)\|_{2}^{2}{\rm d}s{\rm d}t \nonumber\\
\leq&(m_{0}-k_{1})\int_{t_{1}}^{t_{2}}\|\nabla u(t)\|_{2}^{2}{\rm
d}t+(m_{0}-k_{2})\int_{t_{1}}^{t_{2}}\|\nabla
v(t)\|_{2}^{2}{\rm d}t \nonumber\\
\leq&(m_{0}-k)\int_{t_{1}}^{t_{2}}(\|\nabla u(t)\|_{2}^{2}+\|\nabla
v(t)\|_{2}^{2}){\rm d}t  \nonumber\\
\leq&\frac{2(m_{0}-k)(p+2)}{k(p+1)}\int_{t_{1}}^{t_{2}}E(t){\rm d}t.
\label{3.29}
\end{align}
By virtue of \eqref{3.9}, we also have
\begin{align}
&\int_{t_{1}}^{t_{2}}\int_{0}^{t}g_{1}(t-s)\|\nabla
u(t)\|_{2}^{2}{\rm d}s{\rm
d}t+\int_{t_{1}}^{t_{2}}\int_{0}^{t}g_{2}(t-s)\|\nabla
v(t)\|_{2}^{2}{\rm d}s{\rm d}t \nonumber\\
\leq&(m_{0}-k_{1})\int_{t_{1}}^{t_{2}}\|\nabla u(t)\|_{2}^{2}{\rm
d}t+(m_{0}-k_{2})\int_{t_{1}}^{t_{2}}\|\nabla
v(t)\|_{2}^{2}{\rm d}t \nonumber\\
\leq&(m_{0}-k)\int_{t_{1}}^{t_{2}}(\|\nabla u(t)\|_{2}^{2}+\|\nabla
v(t)\|_{2}^{2}){\rm d}t  \nonumber\\
\leq&\frac{2(m_{0}-k)(p+2)}{k(p+1)}\int_{t_{1}}^{t_{2}}E(t){\rm d}t.
\label{3.30}
\end{align}
From \eqref{3.29} and \eqref{3.30}, we see that
\begin{align}
I_{5}=&\frac{1}{2}\int_{t_{1}}^{t_{2}}\int_{0}^{t}g_{1}(t-s)\|\nabla
u(t)-\nabla u(s)\|_{2}^{2}{\rm d}s{\rm
d}t+\frac{1}{2}\int_{t_{1}}^{t_{2}}\int_{0}^{t}g_{2}(t-s)\|\nabla
v(t)-\nabla v(s)\|_{2}^{2}{\rm d}s{\rm d}t \nonumber\\
\leq&\frac{1+2\epsilon}{2}\int_{t_{1}}^{t_{2}}\int_{0}^{t}g_{1}(t-s)\|\nabla
u(t)\|_{2}^{2}{\rm d}s{\rm
d}t+\frac{1}{2}\left(1+\frac{1}{2\epsilon}\right)\int_{t_{1}}^{t_{2}}\int_{0}^{t}g_{1}(t-s)\|\nabla
u(s)\|_{2}^{2}{\rm d}s{\rm d}t \nonumber\\
&+
\frac{1+2\epsilon}{2}\int_{t_{1}}^{t_{2}}\int_{0}^{t}g_{2}(t-s)\|\nabla
v(t)\|_{2}^{2}{\rm d}s{\rm
d}t+\frac{1}{2}\left(1+\frac{1}{2\epsilon}\right)\int_{t_{1}}^{t_{2}}\int_{0}^{t}g_{2}(t-s)\|\nabla
v(s)\|_{2}^{2}{\rm d}s{\rm d}t \nonumber\\
\leq&\frac{1+2\epsilon}{2}\left[\int_{t_{1}}^{t_{2}}\int_{0}^{t}g_{1}(t-s)\|\nabla
u(t)\|_{2}^{2}{\rm d}s{\rm
d}t+\int_{t_{1}}^{t_{2}}\int_{0}^{t}g_{2}(t-s)\|\nabla
v(t)\|_{2}^{2}{\rm d}s{\rm d}t\right] \nonumber\\
&+\frac{1}{2}\left(1+\frac{1}{2\epsilon}\right)\left[\int_{t_{1}}^{t_{2}}\int_{0}^{t}g_{1}(t-s)\|\nabla
u(s)\|_{2}^{2}{\rm d}s{\rm
d}t+\int_{t_{1}}^{t_{2}}\int_{0}^{t}g_{2}(t-s)\|\nabla
v(s)\|_{2}^{2}{\rm d}s{\rm d}t\right] \nonumber\\
\leq&(1+2\epsilon)\frac{(m_{0}-k)(p+2)}{k(p+1)}\int_{t_{1}}^{t_{2}}E(t){\rm
d}t+\left(1+\frac{1}{2\epsilon}\right)\frac{(m_{0}-k)(p+2)}{k(p+1)}\int_{t_{1}}^{t_{2}}E(t){\rm
d}t \nonumber\\
=&\left(2+2\epsilon+\frac{1}{2\epsilon}\right)\frac{(m_{0}-k)(p+2)}{k(p+1)}\int_{t_{1}}^{t_{2}}E(t){\rm
d}t, \label{3.31}
\end{align}
where we have used the following Young's inequalities
\begin{align*}
\left|2\int_{\Omega}\nabla u(t)\nabla u(s){\rm
d}x\right|\leq2\epsilon\|\nabla
u(t)\|_{2}^{2}+\frac{1}{2\epsilon}\|\nabla u(s)\|_{2}^{2},\quad
\forall\ \epsilon>0,
\end{align*}
and
\begin{align*}
\left|2\int_{\Omega}\nabla v(t)\nabla v(s){\rm
d}x\right|\leq2\epsilon\|\nabla
v(t)\|_{2}^{2}+\frac{1}{2\epsilon}\|\nabla v(s)\|_{2}^{2},\quad
\forall\ \epsilon>0.
\end{align*}
 Combining
\eqref{3.23}--\eqref{3.31}, we obtain
\begin{align}
&2\int_{t_{1}}^{t_{2}}E(t){\rm d}t-2(p+1)\int_{t_{1}}^{t_{2}}\int_{\Omega}F(u,v){\rm d}x{\rm d}t \nonumber\\
\leq&2C_{1}E(t_{1})+2C_{p}^{2}E(t_{1})+\left(2+2\epsilon+\frac{1}{2\epsilon}\right)\frac{(m_{0}-k)(p+2)}{k(p+1)}\int_{t_{1}}^{t_{2}}E(t){\rm d}t \nonumber\\
&+C_{3}E(t_{1})+\frac{(m_{0}-k)(p+2)}{k(p+1)}\int_{t_{1}}^{t_{2}}E(t){\rm d}t \nonumber\\
=&C_{0}E(t_{1})+\left(3+2\epsilon+\frac{1}{2\epsilon}\right)\frac{(m_{0}-k)(p+2)}{k(p+1)}\int_{t_{1}}^{t_{2}}E(t){\rm
d}t, \label{3.32}
\end{align}
where $C_{0}=2C_{1}+2C_{p}^{2}+C_{3}.$

On the other hand, by \eqref{3.1} and \eqref{3.9}, we have
\begin{align*}
\int_{t_{1}}^{t_{2}}2(p+1)\int_{\Omega}F(u,v){\rm d}x{\rm
d}t\leq&\int_{t_{1}}^{t_{2}}\frac{p+1}{p+2}\eta(k_{1}\|\nabla
u\|_{2}^{2}+k_{2}\|\nabla v\|_{2}^{2})^{p+2}{\rm d}t \\
\leq&\int_{t_{1}}^{t_{2}}\frac{p+1}{p+2}\eta\left(\frac{2(p+2)}{p+1}E(t)\right)^{p+2}{\rm
d}t
\\
\leq&\int_{t_{1}}^{t_{2}}2\eta\left(\frac{2(p+2)}{p+1}E(0)\right)^{p+1}E(t){\rm
d}t,
\end{align*}
which implies
\begin{align}
&2\int_{t_{1}}^{t_{2}}E(t){\rm d}t-2(p+1)\int_{t_{1}}^{t_{2}}\int_{\Omega}F(u,v){\rm d}x{\rm d}t \nonumber\\
\geq&2\int_{t_{1}}^{t_{2}}E(t){\rm d}t-2\eta\left(\frac{2(p+2)}{p+1}E(0)\right)^{p+1}\int_{t_{1}}^{t_{2}}E(t){\rm d}t \nonumber\\
=&2\left[1-\eta\left(\frac{2(p+2)}{p+1}E(0)\right)^{p+1}\right]\int_{t_{1}}^{t_{2}}E(t){\rm
d}t. \label{3.33}
\end{align}
Note that $E(0)<E_{1},$ we observe that
$$1-\eta\left(\frac{2(p+2)}{p+1}E(0)\right)^{p+1}>0.$$
 Thus, combining \eqref{3.32}
with \eqref{3.33} yields
\begin{align}
2\left[1-\eta\left(\frac{2(p+2)}{p+1}E(0)\right)^{p+1}\right]\int_{t_{1}}^{t_{2}}E(t){\rm
d}t\leq C_{0}E(t_{1})
+\left(3+2\epsilon+\frac{1}{2\epsilon}\right)\frac{(m_{0}-k)(p+2)}{k(p+1)}\int_{t_{1}}^{t_{2}}E(t){\rm
d}t. \label{3.34}
\end{align}
Taking $\epsilon=\frac{1}{2}$ in \eqref{3.34}, we have
\begin{align}
\left[2-2\eta\left(\frac{2(p+2)}{p+1}E(0)\right)^{p+1}-\frac{5(m_{0}-k)(p+2)}{k(p+1)}\right]\int_{t_{1}}^{t_{2}}E(t){\rm
d}t\leq C_{0}E(t_{1}). \label{3.35}
\end{align}

Denote
\begin{align}
\lambda=1- \eta\left(\frac{2(p+2)}{p+1}E(0)\right)^{p+1}-\frac{5(m_{0}-k)(p+2)}{2k(p+1)}.
\label{3.36}
\end{align}
We rewrite \eqref{3.35} as
\begin{align}
2\lambda\int_{t}^{\infty}E(\tau){\rm d}\tau\leq C_{0}E(t),
\label{3.37}
\end{align}
for every $t\in[0,\infty).$

Since $\lambda>0$ when $m_{0}-k>0$ is small enough, by Lemma
\ref{le2.2}, we obtain the following energy decay
$$E(t)\leq E(0)e^{1-2\lambda C_{0}^{-1}t}$$
for every $t\geq C_{0}(2\lambda)^{-1}.$ The proof is completed. \quad
$\Box$

\subsection*{Acknowledgments}
This work was partly supported by the Tianyuan Fund of Mathematics (Grant No. 11026211) and the Natural
Science Foundation of the Jiangsu Higher Education Institutions (Grant No. 09KJB110005).

\end{document}